\documentclass[12pt,reqno,a4paper]{amsart}
\usepackage{amssymb}
\input amssym.def
\usepackage{amsmath,amsfonts,xcolor,mathtools}
\usepackage[bookmarks=false]{hyperref}
\usepackage{amscd}
\usepackage[mathscr]{eucal}
\usepackage{enumitem}
\usepackage{orcidlink}
\allowdisplaybreaks
\numberwithin{equation}{section}

\hypersetup{colorlinks=true,citecolor={purple},linkcolor={teal},urlcolor={violet}}

\theoremstyle{plain}
\newtheorem{theorem}{Theorem}[section]

\newtheorem{lemma}[theorem]{Lemma}

\theoremstyle{definition}

\makeatletter
\@namedef{subjclassname@2020}{%
\textup{2020} Mathematics Subject Classification}
\makeatother



\begin{document}
\title[Congruences for three restricted partition functions of Pushpa and Vasuki]{Congruences modulo powers of $2$ for three restricted partition functions of Pushpa and Vasuki}
\author[Russelle Guadalupe]{Russelle Guadalupe\orcidlink{0009-0001-8974-4502}}
\address{Institute of Mathematics, University of the Philippines Diliman\\
Quezon City 1101, Philippines}
\email{rguadalupe@math.upd.edu.ph}

\begin{abstract}
We establish infinite families of congruences modulo arbitrary powers of $2$ for the three restricted partition functions $M(n), T^\ast(n)$, and $P^\ast(n)$ introduced by Pushpa and Vasuki by employing elementary $q$-series techniques. These generalize some particular congruences for $M(n), T^\ast(n)$, and $P^\ast(n)$ recently found by Nath and Saikia.
\end{abstract}

\maketitle

\section{Introduction}\label{sec1} 

We let $f_m := \prod_{n=1}^\infty(1-q^{mn})$ for a positive integer $m$ and $q\in\mathbb{C}$ with $|q| < 1$ throughout this paper. A partition $\lambda$ of a positive integer $n$ is a nonincreasing sequence of positive integers, known as the parts of $\lambda$, whose sum is $n$. Euler \cite{euler} showed that the generating function for the number $p(n)$ of partitions of $n$ with $p(0) := 1$ is
\begin{align*}
\sum_{n= 0}^\infty p(n)q^n = \dfrac{1}{f_1}.
\end{align*}
Ramanujan \cite{rama} discovered the following remarkable congruences 
\begin{align*}
p(5n+4) &\equiv 0\pmod{5},\\
p(7n+5) &\equiv 0\pmod{7},\\
p(11n+6) &\equiv 0\pmod{11}
\end{align*}
for all $n\geq 0$ using classical $q$-series manipulations. Subsequently, Ramanujan \cite{berno}, Watson \cite{wat}, and Atkin \cite{atk} established congruences modulo powers of $5,7$, and $11$ for $p(n)$.

Pushpa and Vasuki \cite{pushv} introduced in 2022 three restricted partition functions $M(n), T^\ast(n)$, and $P^\ast(n)$ defined as follows:
\begin{enumerate}
\item $M(n)$ is the difference $m_{(d,e)}(n)-m_{(d,o)}(n)$, where $m_{(d,e)}(n)$ (respectively, $m_{(d,o)}(n)$) is the number of partitions of $n$ such that
\begin{itemize}
	\item parts that are congruent to $\pm 1,\pm 3\pmod{10}$ have one color,
	\item parts that are congruent to $0\pmod{10}$ are distinct and have eight colors,
	\item parts that are not congruent to $0,\pm 1,\pm 3\pmod{10}$ are distinct and have four colors each, and 
	\item the number of parts that are not congruent to $\pm 1,\pm 3\pmod{10}$ are even (respectively, odd).
\end{itemize}
\item $T^\ast(n)$ is the difference $t_{(d,e)}(n)-t_{(d,o)}(n)$, where $t_{(d,e)}(n)$ (respectively, $t_{(d,o)}(n)$) is the number of distinct partitions of $n$ with even (respectively, odd) number of parts such that
\begin{itemize}
	\item parts that are congruent to $0\pmod{10}$ have eight colors,
	\item parts that are congruent to $\pm 1,\pm 3\pmod{10}$ have five colors each, and 
	\item parts that are congruent to $\pm 2,\pm 4, 5\pmod{10}$ have four colors each.	
\end{itemize}
\item $P^\ast(n)$ is the difference $p_{(d,e)}(n)-p_{(d,o)}(n)$, where $p_{(d,e)}(n)$ (respectively, $p_{(d,o)}(n)$) is the number of distinct partitions of $n$ with even (respectively, odd) number of parts such that
\begin{itemize}
	\item parts that are congruent to $0\pmod{5}$ have eight colors, and 
	\item parts that are not congruent to $0 \pmod{5}$ have four colors each.	
\end{itemize}
\end{enumerate}

Pushpa and Vasuki \cite{pushv} obtained the generating functions for $M(n), T^\ast(n)$, and $P^\ast(n)$ given by
\begin{align*}
\sum_{n= 0}^\infty M(n)q^n &:= \dfrac{f_2^5f_5^5}{f_1f_{10}},\\
\sum_{n= 0}^\infty T^\ast(n)q^n &:= \dfrac{f_1^5f_{10}^5}{f_2f_5},\\
\sum_{n= 0}^\infty P^\ast(n)q^n &:= f_1^4f_5^4,
\end{align*}
and proved some congruences modulo $2$ and $5$ for $M(n), T^\ast(n)$, and $P^\ast(n)$. Subsequently, Dasappa et. al. \cite{dasa} found more congruences and linear identities for these restricted partition functions. Very recently, Nath and Saikia \cite{naths} established infinite families of congruences modulo powers of $5$ for $M(n)$ and $T^\ast(n)$ and particular congruences for $P^\ast(n)$. They also settled the conjecture of Dasappa et. al. \cite[Conjecture 1]{dasa} by relying on the \textit{Mathematica} package {\tt RaduRK} created by Smoot \cite{smoot}, which is based on Radu's Ramanujan--Kolberg algorithm \cite{radu}. 

We explore in this note more arithmetic properties of $M(n), T^\ast(n)$, and $P^\ast(n)$ using elementary $q$-series techniques. Our main results reveal infinite families of congruences modulo powers of $2$ for these restricted partition functions, which generalizes some of the results of Nath and Saikia \cite{naths}.

\begin{theorem}\label{thm11}
For all integers $k\geq 1$ and $n\geq -1$, we have 
\begin{align}
	M(2^kn+2^{k+1}-1)&\equiv 0\pmod{2^{k-1}},\label{eq11}\\
	T^\ast(2^kn+2^{k+1}-2)&\equiv 0\pmod{2^{k-1}}.\label{eq12}
\end{align}
\end{theorem}

\begin{theorem}\label{thm12}
For all integers $k\geq 0$ and $n\geq -1$, we have
\begin{align}
	P^\ast(2^{4k}n+2^{4k+1}-1)&\equiv 0\pmod{2^{6k}},\label{eq13}\\
	P^\ast(2^{4k+1}n+2^{4k+2}-1)&\equiv 0\pmod{2^{6k+2}},\label{eq14}\\
	P^\ast(2^{4k+2}n+2^{4k+3}-1)&\equiv 0\pmod{2^{6k+3}},\label{eq15}\\
	P^\ast(2^{4k+3}n+2^{4k+4}-1)&\equiv 0\pmod{2^{6k+6}},\label{eq16}\\
	P^\ast(2^{4k+4}n+3\cdot 2^{4k+3}-1)&=0.\label{eq17}
\end{align}
\end{theorem}

We organize the rest of the paper as follows. In Section \ref{sec2}, we present some theta function identities needed to prove Theorems \ref{thm11} and \ref{thm12}, including the generating functions for $M(2n+3), T^\ast(2n+2)$, and $P^\ast(2n+3)$. In Section \ref{sec3}, we demonstrate Theorems \ref{thm11} and \ref{thm12} by finding the generating functions for $M(2^kn+2^{k+1}-1), T^\ast(2^kn+2^{k+1}-2)$, and $P^\ast(2^kn+2^{k+1}-1)$ for any integer $k\geq 1$ via induction.

\section{Preliminaries}\label{sec2}

We provide in this section necessary theta function identities to prove Theorems \ref{thm11} and \ref{thm12}.

\begin{lemma}\label{lem21}
We have the identities
\begin{align}
	\dfrac{f_2^5}{f_{10}} &= \dfrac{f_1^5}{f_5}+5q\dfrac{f_1^2f_2f_{10}^3}{f_5^2},\label{eq21}\\
	\dfrac{f_1}{f_5} &= \dfrac{f_2f_8f_{20}^3}{f_4f_{10}^3f_{40}}-q\dfrac{f_4^2f_{40}}{f_8f_{10}^2},\label{eq22}\\
	f_1f_5^3 &= f_2^3f_{10}-q\dfrac{f_2f_8^2f_{20}^6}{f_4^2f_{10}f_{40}^2}+2q^2f_4f_{20}^3-q^3\dfrac{f_4^4f_{10}f_{40}^2}{f_2f_8^2}.\label{eq23}
\end{align}
\end{lemma}

\begin{proof}
From \cite[Theorem 10.4]{coop}, we have 
\begin{align}
	\dfrac{f_2f_5^5}{qf_1f_{10}^5} &= \dfrac{1}{k(q)}-k(q),\label{eq24}\\
	\dfrac{f_2^4f_5^2}{qf_1^2f_{10}^4} &= \dfrac{1}{k(q)}+1-k(q),\label{eq25}\\
	\dfrac{f_1^3f_5}{qf_2f_{10}^3} &= \dfrac{1}{k(q)}-4-k(q).\label{eq26}
\end{align}
where 
\begin{align*}
	k(q) := q\prod_{n=1}^\infty\dfrac{(1-q^{10n-9})(1-q^{10n-8})(1-q^{10n-2})(1-q^{10n-1})}{(1-q^{10n-7})(1-q^{10n-6})(1-q^{10n-4})(1-q^{10n-3})}	
\end{align*}
is the Ramanujan's function. Subtracting (\ref{eq26}) from (\ref{eq25}) leads to
\begin{align}\label{eq27}
	\dfrac{f_2^4f_5^2}{qf_1^2f_{10}^4} - \dfrac{f_1^3f_5}{qf_2f_{10}^3} = 5.
\end{align}
Multiplying both sides of (\ref{eq27}) by $qf_1^2f_2f_{10}^3/f_5^2$, we obtain (\ref{eq21}). Replacing $q$ with $-q$ in \cite[Theorem 2.1]{hirscs} and using the fact that 
\begin{align*}
	\prod_{n=1}^\infty (1-(-q)^n) = \dfrac{f_2^3}{f_1f_4},
\end{align*}
we obtain (\ref{eq22}). We next infer from (\ref{eq24}) and (\ref{eq25}) that
\begin{align}\label{eq28}
	\dfrac{f_2f_5^5}{qf_1f_{10}^5} = \dfrac{f_2^4f_5^2}{qf_1^2f_{10}^4}-1.
\end{align}
Multiplying both sides of (\ref{eq28}) by $qf_1^2f_{10}^5/(f_2f_5^2)$ and using (\ref{eq22}), we see that
\begin{align}
	f_1f_5^3 &= f_2^3f_{10} - q\dfrac{f_1^2f_{10}^5}{f_2f_5^2}= f_2^3f_{10} - q\dfrac{f_{10}^5}{f_2}\left(\dfrac{f_2f_8f_{20}^3}{f_4f_{10}^3f_{40}}-q\dfrac{f_4^2f_{40}}{f_8f_{10}^2}\right)^2.\label{eq29}
\end{align}
Expanding (\ref{eq29}) yields (\ref{eq23}).
\end{proof}

\begin{lemma}\label{lem22}
Let $q^{-3}f_1^4f_5^4 = \sum_{n=-3}^\infty P^\ast(n+3)q^n$. Then
\begin{align*}
	\sum_{n=-1}^\infty P^\ast(2n+3)q^n = -4\dfrac{f_1^4f_5^4}{q}-8f_2^4f_{10}^4.
\end{align*}
\end{lemma}

\begin{proof}
See \cite[Proposition 3.2]{guad}.
\end{proof}

We now apply Lemmas \ref{lem21} and \ref{lem22} to derive the generating functions for $M(2n+3)$ and $T^\ast(2n+2)$.

\begin{lemma}\label{lem23}
We have that 
\begin{align}
	\sum_{n=-1}^\infty M(2n+3)q^n &= \dfrac{f_1^4f_5^4}{q}-8f_2^4f_{10}^4+10f_1f_2f_5^3f_{10}^3,\label{eq210}\\
	\sum_{n=-1}^\infty T^\ast(2n+2)q^n &= \dfrac{f_1^4f_5^4}{q}+10f_1f_2f_5^3f_{10}^3.\label{eq211}
\end{align}
\end{lemma}

\begin{proof}
Using (\ref{eq21}) and (\ref{eq23}), we have 
\begin{align}\label{eq212}
	\sum_{n=-3}^\infty M(n+3)q^n &= \dfrac{f_2^5f_5^5}{q^3f_1f_{10}} = \dfrac{f_5^5}{q^3f_1}\left(\dfrac{f_1^5}{f_5}+5q\dfrac{f_1^2f_2f_{10}^3}{f_5^2}\right)\nonumber\\
	&= \dfrac{f_1^4f_5^4}{q^3}+5\dfrac{f_1f_2f_5^3f_{10}^3}{q^2}\nonumber\\
	&= \dfrac{f_1^4f_5^4}{q^3}+5\dfrac{f_2f_{10}^3}{q^2}\left(f_2^3f_{10}-q\dfrac{f_2f_8^2f_{20}^6}{f_4^2f_{10}f_{40}^2}+2q^2f_4f_{20}^3-q^3\dfrac{f_4^4f_{10}f_{40}^2}{f_2f_8^2}\right).
\end{align}
We consider the terms involving $q^{2n}$ in the expansion of (\ref{eq212}). In view of Lemma \ref{lem22}, we have 
\begin{align*}
	\sum_{n=-1}^\infty M(2n+3)q^n &= -4\dfrac{f_1^4f_5^4}{q}-8f_2^4f_{10}^4+5\dfrac{f_1^4f_5^4}{q}+10f_1f_2f_5^3f_{10}^3,
\end{align*}
which is exactly (\ref{eq210}).

On the other hand, we apply (\ref{eq21}) and (\ref{eq22}) so that
\begin{align}\label{eq213}
	\sum_{n=-2}^\infty T^\ast(n+2)q^n &= \dfrac{f_1^5f_{10}^5}{q^2f_2f_5} = \dfrac{f_{10}^5}{q^2f_2}\left(\dfrac{f_2^5}{f_{10}}-5q\dfrac{f_1^2f_2f_{10}^3}{f_5^2}\right)\nonumber\\
	&= \dfrac{f_2^4f_{10}^4}{q^2}-5\dfrac{f_1^2f_{10}^8}{qf_5^2}= \dfrac{f_2^4f_{10}^4}{q^2}-5\dfrac{f_{10}^8}{q}\left(\dfrac{f_2f_8f_{20}^3}{f_4f_{10}^3f_{40}}-q\dfrac{f_4^2f_{40}}{f_8f_{10}^2}\right)^2.
\end{align}
Extracting the terms involving $q^{2n}$ in the expansion of (\ref{eq213}) yields (\ref{eq211}).
\end{proof}

\section{Proofs of Theorems \ref{thm11} and \ref{thm12}}\label{sec3}

We establish in this section Theorems \ref{thm11} and \ref{thm12} by using theta function identities and certain generating functions derived from Section \ref{sec2}. As an application of these identities, we derive the following generating functions which are instrumental in the proofs of the aforementioned theorems.

\begin{theorem}\label{thm31}
For all $k\geq 1$, we have that
\begin{align}
	\sum_{n=-1}^\infty M(2^kn+2^{k+1}-1)q^n &= A_k\dfrac{f_1^4f_5^4}{q}-8A_{k-1}f_2^4f_{10}^4+5\cdot 2^kf_1f_2f_5^3f_{10}^3,\label{eq31}\\
	\sum_{n=-1}^\infty T^\ast(2^kn+2^{k+1}-2)q^n &= B_k\dfrac{f_1^4f_5^4}{q}-8B_{k-1}f_2^4f_{10}^4+5\cdot 2^kf_1f_2f_5^3f_{10}^3,\label{eq32}\\
	\sum_{n=-1}^\infty P^\ast(2^kn+2^{k+1}-1)q^n &= C_k\dfrac{f_1^4f_5^4}{q}-8C_{k-1}f_2^4f_{10}^4,\label{eq33}
\end{align}
where for $k\geq 2$, $P_k\in \{A_k, B_k\}$ and $C_k$ satisfy the recurrence 
\begin{align*}
	P_k &= -4P_{k-1}-8P_{k-2}+5\cdot 2^{k-1},\\
	C_k &= -4C_{k-1}-8C_{k-2},
\end{align*}
with initial values $B_0=0, A_0 =A_1=B_1=C_0=1$, and $C_1=-4$.
\end{theorem}

\begin{proof}
We only prove (\ref{eq31}), as the proofs of (\ref{eq32}) and (\ref{eq33}) are similar to that of (\ref{eq31}). We proceed by induction on $k$. In view of Lemmas \ref{lem22} and \ref{lem23}, we note that (\ref{eq31})--(\ref{eq33}) hold for $k=1$. Suppose (\ref{eq31}) holds for some $k\geq 2$. Dividing both sides of (\ref{eq31}) by $q^2$ and employing (\ref{eq23}), we find that
\begin{align}\label{eq34}
	\sum_{n=-3}^\infty &M(2^kn+2^{k+2}-1)q^n \nonumber\\
	&= A_k\dfrac{f_1^4f_5^4}{q^3}-8A_{k-1}\dfrac{f_2^4f_{10}^4}{q^2}+5\cdot 2^k\dfrac{f_1f_2f_5^3f_{10}^3}{q^2}\nonumber\\
	&= A_k\dfrac{f_1^4f_5^4}{q^3}-8A_{k-1}\dfrac{f_2^4f_{10}^4}{q^2}+5\cdot 2^k\dfrac{f_2f_{10}^3}{q^2}\left(f_2^3f_{10}-q\dfrac{f_2f_8^2f_{20}^6}{f_4^2f_{10}f_{40}^2}+2q^2f_4f_{20}^3-q^3\dfrac{f_4^4f_{10}f_{40}^2}{f_2f_8^2}\right).
\end{align}
We read the terms involving $q^{2n}$ in the expansion of (\ref{eq34}). We see from Lemma \ref{lem22} that
\begin{align*}
	\sum_{n=-1}^\infty &M(2^{k+1}n+2^{k+2}-1)q^n \\
	&= A_k\left(-4\dfrac{f_1^4f_5^4}{q}-8f_2^4f_{10}^4\right)-8A_{k-1}\dfrac{f_1^4f_5^4}{q}+5\cdot 2^k\dfrac{f_1^4f_5^4}{q}+5\cdot 2^{k+1}f_1f_2f_5^3f_{10}^3\\
	&= (-4A_k-8A_{k-1}+5\cdot 2^k)\dfrac{f_1^4f_5^4}{q}-8A_kf_2^4f_{10}^4+5\cdot 2^{k+1}f_1f_2f_5^3f_{10}^3\\
	&= A_{k+1}\dfrac{f_1^4f_5^4}{q}-8A_kf_2^4f_{10}^4+5\cdot 2^{k+1}f_1f_2f_5^3f_{10}^3,
\end{align*}
where the last equality follows from the definition of $A_k$. Thus, (\ref{eq31}) holds for $k+1$, so (\ref{eq31}) holds for all $k\geq 1$ by induction.
\end{proof}

\begin{proof}[Proof of Theorem \ref{thm11}]
We claim that for all $k\geq 1$ and $P_k\in \{A_k, B_k\}$, $P_k$ is divisible by $2^{k-1}$ but not by $2^k$. We prove this via strong induction on $k$. We compute $A_2=6$ and $B_2=-2$, so the claim holds for $k\in \{1,2\}$. Suppose now that the claim holds for all positive integers at most some $k\geq 3$. Write $P_k = 2^{k-1}a_k$ and $P_{k-1}=2^{k-2}b_k$ for some odd integers $a_k$ and $b_k$. We infer from the definition of $P_k$ that 
\begin{align*}
	P_{k+1} = -4P_k-8P_{k-1}+5\cdot 2^k = 2^k(-2a_k-2b_k+5).
\end{align*}
Since $5-2a_k-2b_k$ is odd, we obtain $2^k\mid P_k$ but $2^{k+1}\nmid P_k$. Thus, the claim holds for $k+1$, and hence for $k\geq 1$ by strong induction. Combining the above claim and the generating functions (\ref{eq31}) and (\ref{eq32}), we arrive at (\ref{eq11}) and (\ref{eq12}), respectively.
\end{proof}

To prove Theorem \ref{thm12}, we require the following result on the exact values of $C_k$. 

\begin{lemma}\label{lem32}
For all integers $k\geq 0$, we have 
\begin{align*}
	C_{4k} &= (-64)^k,\\
	C_{4k+1} &= -4(-64)^k,\\
	C_{4k+2} &= 8(-64)^k,\\
	C_{4k+3} &= 0
\end{align*}
\end{lemma}

\begin{proof}
Using the definition of $C_k$, we observe that for all $k\geq 0$ and $r\in \{0,1,2,3\}$,
\begin{align}
	C_{4k+r+4}&+64C_{4k+r}\nonumber\\
	&= C_{4k+r+4}+4C_{4k+r+3}+8C_{4k+r+2}-4(C_{4k+r+3}+4C_{4k+r+2}+8C_{4k+r+1})\nonumber\\
	&+8(C_{4k+r+2}+4C_{4k+r+1}+8C_{4k+r})= 0.\label{eq35}
\end{align}
We see from (\ref{eq35}) that $C_{4k+r} = (-64)^kC_r$ for all $k\geq 0$ and $r\in \{0,1,2,3\}$. In view of the computed values $C_2=8$ and $C_3=0$, we arrive at the desired values of $C_k$. 
\end{proof}

\begin{proof}[Proof of Theorem \ref{thm12}]
Appealing to Lemma \ref{lem32}, we first replace $k$ with $4k$ in (\ref{eq33}) so that
\begin{align*}
	\sum_{n=-1}^\infty P^\ast(2^{4k}n+2^{4k+1}-1)q^n &= (-64)^k\dfrac{f_1^4f_5^4}{q},
\end{align*}
yielding (\ref{eq13}).
We next replace $k$ with $4k+1$ in (\ref{eq33}) so that
\begin{align*}
	\sum_{n=-1}^\infty P^\ast(2^{4k+1}n+2^{4k+2}-1)q^n &= 4(-64)^k\dfrac{f_1^4f_5^4}{q}-8(-64)^kf_2^4f_{10}^4,
\end{align*}
which gives (\ref{eq14}).
We next replace $k$ with $4k+2$ in (\ref{eq33}) so that
\begin{align*}
	\sum_{n=-1}^\infty P^\ast(2^{4k+2}n+2^{4k+3}-1)q^n &= 8(-64)^k\dfrac{f_1^4f_5^4}{q}+32(-64)^kf_2^4f_{10}^4,
\end{align*}
which gives (\ref{eq15}). 
We now replace $k$ with $4k+3$ in (\ref{eq33}) so that
\begin{align}\label{eq36}
	\sum_{n=-1}^\infty P^\ast(2^{4k+3}n+2^{4k+4}-1)q^n &= (-64)^{k+1}f_2^4f_{10}^4,
\end{align}
arriving at (\ref{eq16}). Dividing both sides of (\ref{eq36}) by $q$, we have
\begin{align}\label{eq37}
	\sum_{n=-2}^\infty P^\ast(2^{4k+3}n+3\cdot 2^{4k+3}-1)q^n &= (-64)^{k+1}\dfrac{f_2^4f_{10}^4}{q}.
\end{align}
Examining the terms involving $q^{2n}$ in the expansion of (\ref{eq37}), we finally obtain (\ref{eq17}).
\end{proof}

%
%
%
%
%
%
%


\begin{thebibliography}{99}

\bibitem{atk} A. O. L. Atkin, Proof of a conjecture of Ramanujan, {\it Glasg. Math. J.} {\bf 8} (1967), 14--32.	

\bibitem{berno} B. C. Berndt and K. Ono, Ramanujan's unpublished manuscript on the partition and tau functions with proofs and commentary, in {\it The Andrews Festschrift} (Springer-Verlag, Berlin, 2001), pp. 39--110.

\bibitem{coop} S. Cooper, {\it Ramanujan's Theta Functions} (Springer, Cham, 2017).

\bibitem{dasa} R. Dasappa, Channabasavayya and G. K. Keerthana, On some new arithmetic properties of certain restricted color partition functions, {\it Arab. J. Math.} {\bf 13} (2024), 275--289.

\bibitem{euler} L. Euler, {\it Introductio in Analysin Infinitorum} (Marcum-Michaelem Bousquet, Lausannae, 1748).

\bibitem{guad} R. Guadalupe, Linear identities for partition pairs with $5$-cores (2026), preprint, arXiv:2601.04743.

\bibitem{hirscs} M. D. Hirschhorn and J. A Sellers, Elementary proofs of parity results for $5$-regular partitions, {\it Bull. Aust. Math. Soc.} {\bf 81} (2010), 58--63.

\bibitem{naths} H. Nath and M. P. Saikia, Arithmetic properties of partition functions introduced by Pushpa and Vasuki, {\it J. Symbolic Comput.} {\bf 135} (2026), paper no. 102555, 15 pp.

\bibitem{pushv} K. Pushpa and K. R. Vasuki, On Eisenstein series, color partition and divisor function, {\it Arab. J. Math.} {\bf 11} (2022), 355--378.

\bibitem{radu} C.-S. Radu, An algorithmic approach to Ramanujan–Kolberg identities, {\it J. Symbolic Comput.} {\bf 68} (2015), 225--253.

\bibitem{rama} S. Ramanujan, Some properties of $p(n)$, the number of partitions of $n$, {\it Proc. Cambridge Math. Soc.} {\bf 19} (1919), 210--213.

\bibitem{smoot} N. A. Smoot, On the computation of identities relating partition numbers in arithmetic progressions with eta quotients: An implementation of Radu's algorithm, {\it J. Symbolic Comput.} {\bf 104} (2021), 276--311.

\bibitem{wat} G. Watson, Ramanujans Vermutung \"{u}ber Zerf\"{a}llungszahlen, {\it J. Reine Angew. Math.} {\bf 179} (1938), 97--128.
\end{thebibliography}

\end{document}